# A DIFFERENT PERSPECTIVE ON TEACHING GEOMETRY AT HIGH SCHOOL THE GREEK CASE STUDY

Rizos, I.[1]
Gkrekas, N.[2]
[1]Department of Mathematics, University of Thessaly, Lamia 35100, Greece
(ioarizos@uth.gr) https://orcid.org/0000-0002-4092-1715
[2]Department of Mathematics, University of Thessaly, Lamia 35100, Greece
(ngkrekas@uth.gr) https://orcid.org/0000-0001-9665-4559

Geometry is essentially a global language, which is fully understood in different times, countries and cultures. The proof of a geometric theorem (e.g. the Pythagorean Theorem) or a geometric construction (e.g. the construction of an equilateral triangle with straightedge and compass) is recognized in more or less the same way by everyone. This fact, in combination with other parameters, brings out a human face for geometric and mathematical achievements in general, thus allowing people and civilizations to come together.

The pedagogical value of teaching Geometry in Secondary Education is indisputable, mainly because it helps students to develop their ability to perceive space and directly links Mathematics with the real world. At the same time, as a school course, it expresses an environment in which the initiation into the proof process and the understanding of how a basic branch of mathematics is founded and developed thrives. Nevertheless, school experience shows that the teaching of Geometry, especially in high school (10th-12th grade), presents significant challenges, with the most characteristic the understanding of proofs [3], [12] and their mechanistic execution [6].

Particularly in Greece, the gradual shrinking of the Euclidean Geometry course and the static approach to the development of its remaining content, weaken the arguments for challenging the extreme position of Jean Dieudonné, who 63 years ago, with the aim of reforming the curricula for Mathematics, proclaimed: "Euclid must go!" [2, p. 35]. In practice, Euclidean geometry is taught detached from its physical meaning and largely unrelated to students' experiences. As a consequence, this gives the impression that it is a closed logical system, like a game with peculiar rules (e.g. the chess).

At the same time, the formalistic character of school mathematics, i.e. the focus on techniques for solving theoretical exercises [9], combined with the pragmatic perception that pervades many textbooks and articles, even the so-called "popularized" ones, do not favor the development of different views on Euclidean geometry. Thus, a formalistic model of teaching and learning is essentially adopted, which privileges teacher-centeredness over discovery and exchanges conceptual understanding for mechanistic task-solving and algorithmic thinking.

On the other hand, the curriculum and the curricular guidelines –especially during the Covid-19 pandemic era and the subsequent implementation of distance education, when the teaching model changed radically (see [7])– impose strict time constraints on teachers, who daily try to strike a balance between covering the syllabus and students' understanding. In such a complex environment, therefore, the possibility of devoting teaching time to the types of geometries (Euclidean, Non-Euclidean, Pseudo-Euclidean etc.) and their methods of foundation may not seem realistic.

However, research has shown [1], [4], [10], [11], [5], [8] that by coming with appropriate resources (e.g. teaching scenarios or thought experiments) in contact with geometries beyond Euclidean and by discovering the meaning and role of metrics, students have the opportunity to challenge entrenched notions about the nature of geometry itself, to interact, negotiate mathematical meanings, discover connections with Physics and make transgressions, thus shaping appropriate conditions and opening horizons for further learning and understanding in Mathematics.

In the above context it may also be possible to see the three basic geometries of Euclid, Galileo and Minkowski (see [13]), and to integrate elements of them into teaching. This may give students and/or undergraduate students the future perspective of studying a) the way of viewing and the role of postulates, b) the fundamental concept of "space" and its modeling, c) the interaction of mathematical viewing with physical theory, d) the importance of geometric transformations and their invariants, hence of different representations of a shape (e.g. the "circle") from geometry to geometry, and e) episodes from the history of mathematics.

The above can help to shift the horizon of children's thinking and, to some extent, change the way they perceive mathematical patterns and the functioning of the universe. Pupils are released from the obligation to accept as axiomatic something that can be derived from the study of physical phenomena. After all, as the History of Science reminds us, the approach to geometry through experiment and observation, combined with the long philosophical and epistemological transformation of the concept of 'space', led to scientific revolutions and the transcendence of established philosophical beliefs about the natural world. It would therefore be interesting, in the future, to conduct a study that would examine whether and in what ways the teaching of the historical evolution of geometric concepts and processes can create appropriate conditions for high school students to acquire positive attitudes towards geometry and to become actively involved in the negotiation of mathematical meanings.